\documentclass{elsart}
\usepackage{epsfig,graphicx,graphics,latexsym,amssymb,amsfonts,amsmath, amscd}
\usepackage{enumerate}
\usepackage{euscript}
\usepackage{float}
\newcommand{\Z}{\hbox{Z\hskip -4pt Z}}

\newtheorem{definition}{{\bf Definition}}[section]
\newtheorem{theorem}[definition]{{\bf Theorem}}

\newtheorem{proposition}[definition]{\noindent {\bf Proposition}}
\newtheorem{lemma}[definition]{\noindent {\bf Lemma}}

\newtheorem{problem}[definition]{\noindent {\bf Problem}}

\def\Proof{{\parindent0pt  {\bf Proof.\ }}}
\def\endproof{\hfill {\kern 6pt\penalty 500
\raise -0pt\hbox{\vrule \vbox to5pt {\hrule width 5pt
\vfill\hrule}\vrule}}}

\footnotetext

\begin{document}
\begin{frontmatter}
\title{Claw-freeness, $3$-homogeneous subsets of a graph and  a reconstruction problem\thanksref{label1}}

\author{Maurice Pouzet},
%\address{ICJ, Universit\'e Claude Bernard Lyon1,
% 43 Bd.  11 Novembre 1918, 69622 Villeurbanne Cedex, France}
\thanks[label1]{Done under the auspices of the French-Tunisian CMCU "Outils math\'ematiques pour l'Informatique" 05S1505}
\ead{pouzet@univ-lyon1.fr}
\author{Hamza  Si Kaddour\corauthref{cor}}
\address{Institut Camille-Jordan, Math\'ematiques, Universit\'e Claude Bernard Lyon1,
 43 Bd.  11 Novembre 1918, 69622 Villeurbanne Cedex, France}
\ead{sikaddour@univ-lyon1.fr}
\corauth[cor]{Corresponding author.}

\date{\today}
\begin{center}Mailbox
\end{center}
\begin{abstract}
We describe
$Forb\{K_{1,3}, \overline {K_{1,3}}\}$, the class of graphs $G$ such that 
 $G$ and its complement 
$ \overline{G}$ are claw-free. With few exceptions, it is made  of graphs whose connected components consist of cycles of length 
at least 4, paths or isolated vertices, and of the complements of these graphs. Considering the hypergraph ${\mathcal H} ^{(3)}(G)$ made of the $3$-element subsets of the vertex set of a graph $G$ on
which $G$ induces  a clique or an independent subset, we deduce  from above a  description of the Boolean sum $G\dot{+}G'$  of two graphs  $G$ and $G'$ giving  the same hypergraph. We indicate the role of  this latter description in a reconstruction problem of graphs up to complementation.
\end{abstract}
%\noindent {\it MSC:} 04A05; 05C60.\\
\end{frontmatter}

\section{Results and motivation} 
  Our notations and terminology mostly follow \cite {Bo}. The graphs we consider in this paper are  undirected, simple and have no  loop.   That is  a {\it graph} is a pair $G:= (V, \mathcal E)$, where $\mathcal E$ is a subset of $[V]^2$, the set of $2$-element subsets of $V$. Elements of $V$ are the {\it vertices} of $G$ and elements of $\mathcal E$ its {\it edges}.  We denote by $V(G)$ the vertex set of $G$ and  by $E(G)$ its edge set. We look at members of $[V]^2$ as unordered pairs of distinct vertices. If $A$ is a subset of $V$, the pair $G_{\restriction A}:=(A, \mathcal E\cap [A]^2)$ is the \emph{graph induced by $G$ on $A$}. The {\it complement} of $G$ is the simple graph 
 ${\overline G}$ whose vertex set is $V$ and whose edges are the unordered pairs of nonadjacent and distinct vertices of $G$, that is $\overline G =(V, {\overline {\mathcal E}})$, where ${\overline {\mathcal E}}=[V]^2\setminus \mathcal E$. %and by $e(G):=\vert E(G)\vert $ the number of edges. If $\{x,y\}$ is an edge of $G$ we set $G(x,y)=1$; otherwise we set $G(x,y)=0$.  For $A\subseteq V$ and $\varepsilon \in \{0,1\}$, $G(x,A)=\varepsilon$ means $G(x,a)=\varepsilon$
%for all $a\in A$. The {\it degree} of a vertex $x$  of $G$, denoted $d_G(x)$, is  the number of  edges which contain $x$. The graph $G$ is {\it regular} if $d_G(x)=d_G(y)$ for all $x,y\in V$.  If $K$ is a subset of $V$, the {\it  restriction} of $G$ to $K$, also called the {\it induced graph} on $K$ is the graph $G_{\restriction K}:= (K, [K]^2\cap \mathcal E)$. 
%The induced graph $G_{\restriction {V\setminus K}}$ is denoted $G_{\restriction {- K}}$, and more simply $G_{ {- x}}$ if $K=\{x\}$.
We denote by 
$K_3$  the complete graph on
$3$ vertices and  by $K_{1,3}$ the graph  made of a
vertex linked to a
 $\overline {K_{3}}$, this graph is called a \emph{claw}. We denote by $A_{6}$  the graph on $6$ vertices made of a $K_3$  bounded by
three $K_3$ (cf. Figure 1) and by $C_n$ the $n$-element cycle, $n\geq 4$. We denote by
$P_9$ the Paley graph on $9$ vertices (cf. Figure 1). Note that $P_9$ is isomorphic to  its complement
$\overline{P_9}$ and also to $K_3\Box K_3$, the  cartesian  product of $K_3$ by
itself  (see \cite{Bo} page 30 if needed for a definition of the \emph{cartesian product of graphs}, and  see \cite{VW} 
page 176 and  \cite {Bo} page 28 for a definition and basic properties of \emph{Paley graphs}).

\begin{figure}[H]
\begin{center}
\includegraphics[width=4.5in]{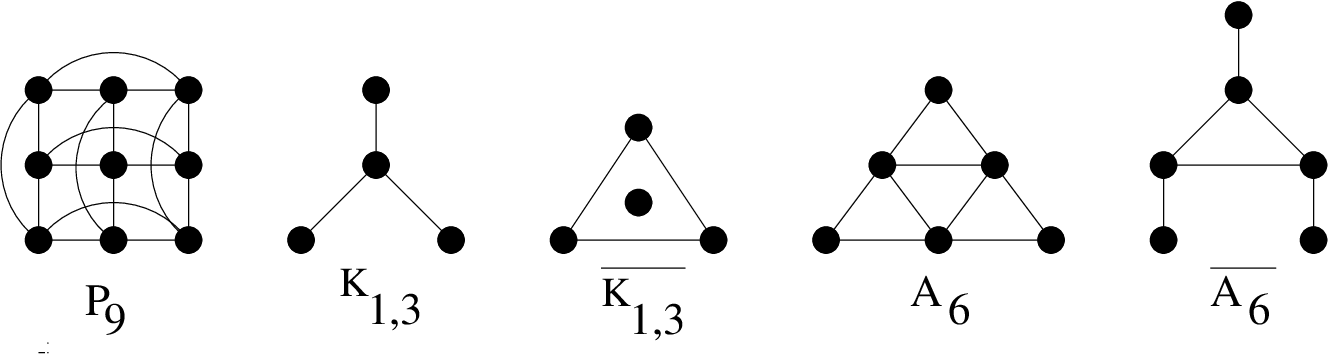}
\end{center}
\caption{} \label{1}
\end{figure}

Given a set  $\mathcal F$ of graphs, we denote by
$Forb \mathcal F$ the class of graphs $G$ such that
no member of $\mathcal F$ is isomorphic to an
induced subgraph of $G$.  Members of
$Forb \{K_{3}\}$, resp. $Forb\{K_{1,3}\}$ are called  {\it triangle-free}, resp. {\it
claw-free}  graphs.
The main result of this note asserts: 

 \begin{theorem} \label{K_{1,3}} The class 
$Forb\{K_{1,3}, \overline {K_{1,3}}\}$ consists of $A_6$, of the induced
subgraphs of $P_9$, of graphs whose connected components consist of 
cycles of length at least $4$, paths or isolated
vertices, and of the complements of these graphs.
\end{theorem}
As an immediate consequence of Theorem \ref {K_{1,3}}, note that the
graphs 
 $A_{6}$ and  $\overline
{A_{6}}$ are the only members of $Forb\{K_{1,3}, \overline
{K_{1,3}}\}$ which contain a $K_{3}$  and  a $\overline
{K_{3}}$ with no vertex in common.

From  Theorem \ref{K_{1,3}} we  obtain a characterization of  the Boolean sum of two graphs having the same $3$-homogeneous subsets. For that, we say that a  subset of vertices of a graph $G$ is {\it
homogeneous} if it is a clique or an independent set (note that if the word homogeneous is used with this meaning in Ramsey  theory;  in other areas of graph theory it has other meanings, several in fact).
Let ${\mathcal H}^{(3)}(G)$ be the hypergraph having the same
vertices as
$G$ and whose  
hyperedges are the 
$3$-element homogeneous subsets of $G$. Given two graphs $G$ and $G'$
on the same vertex set $V$, we recall that the {\it
Boolean sum}
$G\dot{+}G'$ of $G$ and $G'$ is the graph  on $V$ whose edges are unordered pairs $e$ of
distinct vertices such that $e\in E(G)$ if and only if $e\notin E(G')$. Note that $E(G\dot{+}G')$ is the symmetric difference $E(G)\Delta E(G')$ of $E(G)$ and $E(G')$. The graph $G\dot{+}G'$ is also called the \emph{symmetric difference} of $G$ and $G'$ and denoted by $G\Delta G'$ in \cite{Bo}. 

Given a graph $U$ with vertex set  $V$, the {\it
edge-graph} of
$U$ is the graph
$S(U)$ whose
vertices are the edges $u$ of $U$ and whose edges are  unordered pairs $\{u,v\}$
such that
$u=\{x,y\}$, $v=\{x,z\}$ for three distinct elements $x,y,z\in V$ such that
$\{y,z\}$  is
not an edge of $U$. Note that the edge-graph $S(U)$ is a spanning subgraph of $L(U)$, the \emph{line-graph} of $U$, not to be confused with it.

Claw-free graphs and triangle-free graphs are related by means 
of the edge-graph
construction. Indeed, as it is immediate to see, for every graph $U$, we have:
\begin{equation}\label{eqclawfree}
U\in Forb \{K_{1,3}\}\Longleftrightarrow S(U)\in Forb \{K_{3}\}.
\end{equation} 
 
Our characterization is this:

\begin{theorem}\label{S(U)}
Let $U$ be a graph. The following properties are equivalent:
\begin{enumerate}[{(1)}]
\item 
There are two graphs $G$ and $G'$ having the same
$3$-element homogeneous subsets such that $U:=G\dot{+} G'$;
\item $S(U)$ and $S(\overline U)$ are bipartite;
\item
\begin{enumerate}[{(i)}]
\item Either $U$ is  an induced subgraph of $P_9$\label{P9}, 
\item or  the connected components of
$U$,  or of its complement $\overline U$,  are cycles of even length, paths or isolated vertices\label{direct}.
 \end{enumerate}
  \end{enumerate}
\end{theorem}

As a consequence,  if  the  graph  $U$ satisfying Property (1) is disconnected, then $U$ contains no $3$-element cycle\label{claw}, moreover, 
if $U$ contains no $3$-element cycle then each connected component of $U$ is  a cycle of even length,  a path, or an isolated vertex\label{degre2}, in particular  $U$ is bipartite.\label{bipartite}

 The implication $(2)\Rightarrow  (3)$  in Theorem \ref{S(U)} follows immediately from Theorem \ref {K_{1,3}}. Indeed,
suppose  that Property (2) holds, that is  $S(U)$ and $S(\overline U)$ are bipartite, then
from Formula (\ref{eqclawfree}) and from the fact that $S(A_6)$ and 
$S(C_n)$, $n\geq 4$, are respectively isomorphic to $C_9$ and to $C_n$, we have: $$U\in   
Forb\{K_{1,3}, \overline {K_{1,3}}, A_6,\overline {A_6}, C_{2n+1}, \overline {C_{2n+1}}
: n\geq 2\}.$$ From Theorem \ref{K_{1,3}}, Property (3) holds. The other implications, obtained by more straigthforward arguments, are given in Subsection \ref{Proof of TheoremS(U)}.

This leaves open the following:
\begin{problem}
Which pairs of graphs $G$ and $G'$  with the same
$3$-element homogeneous subsets have a given Boolean sum $U
:=G\dot{+} G'$?
\end{problem}

A partial answer, motivated by the reconstruction problem discussed below, is given in \cite{dlps2}. We  mention  that two graphs $G$ and $G'$ as above are determined  by the graphs induced on the  connected components of $U
:=G\dot{+} G'$  and on a system of distinct representatives of these connected components (Proposition 10 \cite{dlps2}).

A quite natural problem, related to the study of Ramsey numbers for triples,  is this:
\begin{problem} Which hypergraphs are of the form ${\mathcal H}^{(3)}(G)$ ?
\end{problem} 
 An asymptotic  lower bound of the  size of ${\mathcal H}^{(3)}(G)$ in terms of $\vert V(G)\vert $ was established by A.W.Goodman \cite{goo}.

The motivation for Theorem 
 \ref{S(U)} (and thus Theorem \ref {K_{1,3}}) originates in a reconstruction problem on graphs that we present now.  
 Considering two graphs $G$ and $G'$ on the same set $V$ of vertices, we say that $G$ and $G'$ are {\it isomorphic up to complementation} if $G'$ is isomorphic to $G$ or to the complement $\overline G$ of $G$. Let $k$ be a non-negative integer,  we say that $G$ and $G'$ are {\it $k$-hypomorphic  up to complementation} if for every $k$-element subset $K$ of $V$, the graphs  $G_{\restriction K}$ and $G'_{\restriction K}$  induced by $G$ and $G'$ on $K$ are isomorphic up to complementation. Finally, we say that $G$ is {\it $k$-reconstructible up to complementation} if every graph  $G'$ which is $k$-hypomorphic to   $G$ up to complementation is in fact isomorphic to $G$ up to complementation. 
 The following problem emerged from a question of P.Ille \cite{Ille}:
 \begin{problem} For which pairs $(k,v)$ of integers, $k<v$,   every graph  $G$  on $v$ vertices is $k$-reconstructible up to complementation?
\end{problem}

It is immediate to see that  if the conclusion of the problem above is positive,  $v$ is distinct from $3$ and $4$ and, with a little bit of thought, that if $v\geq 5$ then $k\geq 4$ (see Proposition 4.1 of \cite{dlps1}). 
 With  J.
Dammak, G. Lopez \cite{dlps1} and  \cite{dlps2} we proved that the conclusion is positive if: \begin {enumerate} [{(i)}]
\item $4\leq k\leq v-3$  or 
\item $4\leq k=v-2$ and $v\equiv 2\;  (mod \; 4)$.  \end{enumerate}
We do not know if in (ii) the condition $v\equiv 2\;  (mod \; 4 )$ can be dropped. For  $4\leq k=v-1$, we checked that the conclusion holds if $v=5$ and noticed that for larger values of $v$ it could be
negative or extremely hard to obtain, indeed, a positive conclusion  would imply that Ulam's reconstruction
conjecture holds (see Proposition 19 of \cite{dlps2}).

The reason for which Theorem \ref{S(U)} plays a role in that matter relies on the following result, an easy consequence of the famous  Gottlieb-Kantor theorem on incidence matrices (\cite{Go,KA}). 

\begin{proposition} \label{down}(Proposition 2.4 \cite{dlps1}) Let  $t \leq min{(k,  v-k)}$ and $G$ and $G'$ be two graphs on the same set $V$ of $v$ vertices.  If $G$ and $G'$ are $k$-hypomorphic up to complementation then they are $t$-hypomorphic up to complementation.
\end{proposition}

Indeed, if $3\leq k\leq v-3$, Proposition \ref{down} tells us that two graphs $G$ and $G'$ which  are $k$-hypomorphic up to complementation are $3$-hypomorphic up to complementation, which amounts to the fact that $G$ and $G'$ have the same $3$-homogeneous subsets. A carefull study of such pairs $G$, $G'$ allows to deal with the case $v=k+3$ and $v\equiv 1\; mod\;4$ (see \cite{dlps2}). Other cases use
properties of the rank of some incidence  matrices; notably  a result of  R.M. Wilson 
\cite{W} on incidence matrices (and also the Gottlieb-Kantor theorem). In these cases the conclusion is stronger: $G'$ or $\overline{G'}$  is equal
to $G$ (see \cite{dlps1}). 
\section{Proofs}\label{proofs}
Let $U$ be a graph.
For an unordered pair $e:=\{x,y\}$ of distinct vertices, we set $U(e)=1$ if
$e\in E(U)$ and $U(e)=0$ otherwise. 
Let   $x\in V(U)$; we
denote by 
$N_U(x)$ and $d_{U}(x)$ the {\it neighborhood} and
the {\it degree} of
$x$ (that is $N_U(x):=\{
y\in V(U) :
\{x, y\}\in E(U)\}$ and $d_U(x):=\vert N_U(x)\vert$). 
Set $U_{(2)}:= \{ x \in V(U)  : d_{ U}(x)\leq
2
\}$, thus
$\overline {U}_{(2)}= \{ x \in V(U)  : d_{\overline
U}(x)\leq 2\}$. 

\subsection {Proof of Theorem \ref{K_{1,3}}.}
Trivially, the graphs described in Theorem \ref{K_{1,3}} belong to $Forb\{K_{1,3}, \overline {K_{1,3}}\}$.
We deduce the converse  from Lemmas \ref {trivial} and  \ref {caseanalysis} below.

 \begin{lemma}\label {trivial}
Let $U$ be  graph,  $U\in Forb\{K_{1,3}, K_{3}\}$ if and only if  the connected components of $U$ are   cycles of length at least $4$, paths or
isolated vertices. 
\end{lemma}
 The \emph{bull} is the graph $B_5$ on $5$ vertices made of  a $K_3$ and two additional edges linked to  $K_3$ and with no common vertex;  let $E_6$ be the graph on $6$ vertices made by  a square
 bounded by two $K_3$ with  a common vertex. If $v$, $w$ are two vertices of $P_{9}$, we denote by $P_{9-v}$, resp. $P_{9-vw}$,   the graph obtained from $P_{9}$ by deleting the vertex $v$, resp. the vertices   $v$ and $w$. Since the automorphisms of $P_9$ act transitively on the vertices, on the edges and on the non-edges, the graph $P_{9-v}$, considered up to an isomorphism, does not depend upon the vertex $v$, similarly, if $v$ and $w$ are  distinct, the graph $P_{9-vw}$ only depends upon the fact that the unordered pair $\{u, v\}$ is or is not an edge. In order to distinguish between these two cases, we will denote by $P_{9-\varepsilon}$ the graph obtained by deleting the edge $\varepsilon$.  Note that $B_5$ and $P_{9-v}$ are isomorphic to their  complements and also note that $B_5$, $E_6$, $\overline {E_{6}}$, $P_{9-\varepsilon}$, $\overline {P_{9-\varepsilon}}$,  $P_{9-v}$  are induced subgraphs
of $P_{9}$.These graphs are represented Figure 2.

 \begin{figure}[H]
 \begin{center}
\includegraphics[width=4.7in]{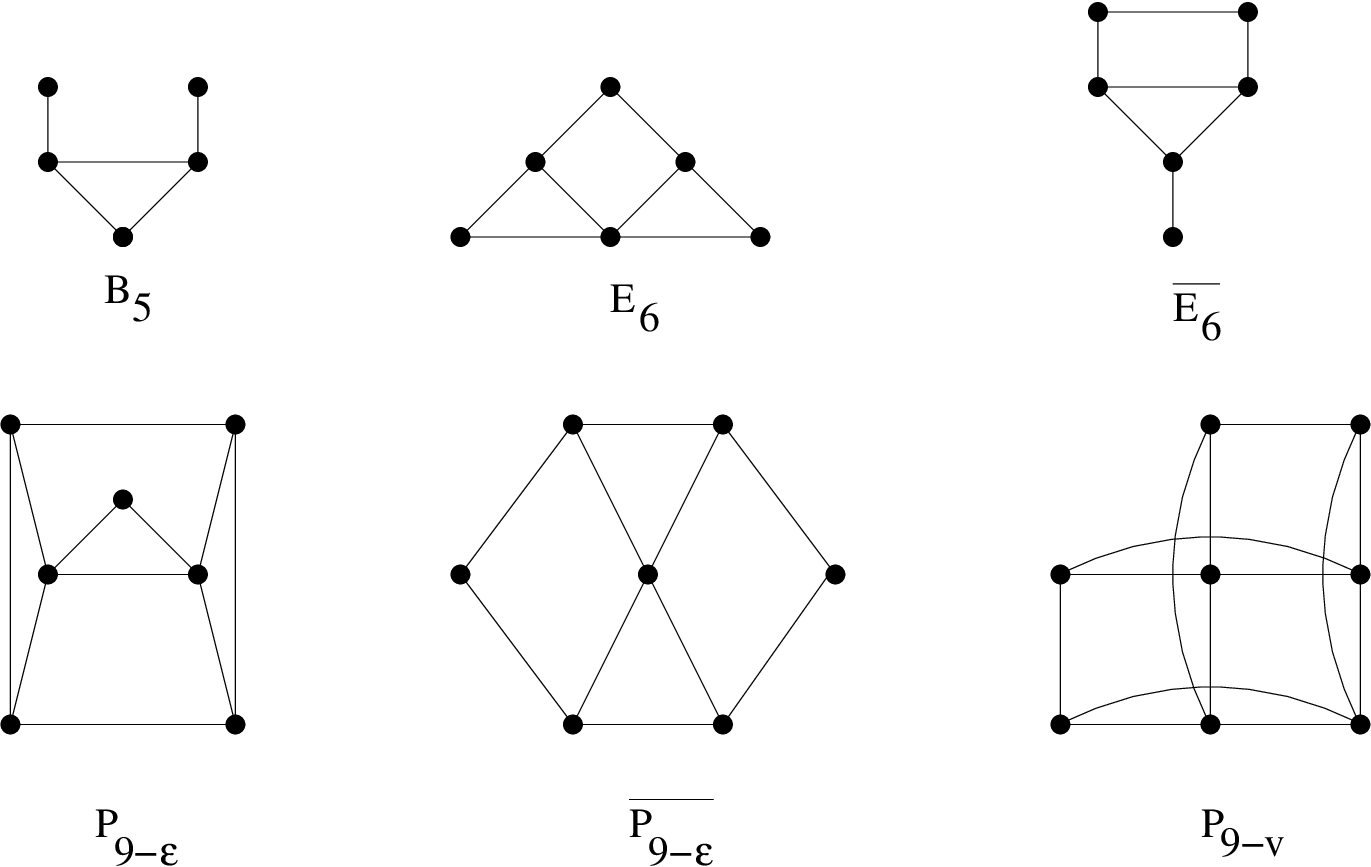}
\caption{} \label{2}
\end{center}
\end{figure}
  \begin{lemma} \label{caseanalysis} Let $U\in
Forb\{K_{1,3}, {\overline {K_{1,3}}}\}$. If $U$ contains
a
$K_{3}$  and a
$\overline {K_{3}}$ then  $U$ is isomorphic to one of the
following
 nine graphs:
the bull $B_{5}$, $A_{6}$,  $\overline {A_{6}}$, $E_{6}$, $\overline {E_{6}}$, $P_{9-\varepsilon}$, 
 $\overline {P_{9-\varepsilon}}$, $P_{9-v}$, $P_{9}$.
\end{lemma}
With these two Lemmas, the proof of Theorem \ref{K_{1,3}} is immediate.  Indeed,  let   $U\in
Forb\{K_{1,3}, {\overline {K_{1,3}}}\}$. If $U$ contains no $K_3$ then according to Lemma \ref{trivial}, the connected components of $U$ are  cycles of length at least $4$, paths or
isolated vertices,  thus $U$ has the form announced  in Theorem  \ref{K_{1,3}}. If $U$ contains no $\overline {K_{3}}$ then via  Lemma \ref{trivial} applied to $\overline U$ we get the same conclusion.   If $U$ contains
a
$K_{3}$  and a
$\overline {K_{3}}$ then  Lemma \ref{caseanalysis} tells us that except
$A_{6}$ and  $\overline {A_{6}}$, each of the graph listed above  is
isomorphic  to an induced subgraph of  $P_{9}$.  Thus $U$ has the form announced  in Theorem  \ref{K_{1,3}}. 

\subsubsection{Proof of Lemma  \ref{trivial}.} The proof is a consequence of the following trivial observation.  
\begin {claim} \label {triangle}  If  $U\in Forb\{K_{1,3}\}$, every $x\in V(U)\setminus  U_{(2)}$
belongs to a $K_{3}$.
\end{claim}
Indeed,  let $U\in Forb\{K_{1,3}, K_{3}\}$. According to Claim 
\ref{triangle},  $V(U)= U_{(2)}$, thus the connected components of $U$ are cycles, paths or
isolated vertices, and  since $U\in Forb\{ K_{3}\}$, the cycles have length at least four. The converse is immediate. 
 
\subsubsection{ Ingredients of the proof of Lemma
\ref{caseanalysis}.}
 
 The proof of Lemma
\ref{caseanalysis} given in Subsection \ref{Prooflemmacaseanalysis} below is a case by case analysis using the
following claims. 
\begin{claim} \label {disjointtriangle} An arbitrary graph $U$ which
contains a
$K_{3}$ and a 
$\overline {K_{3}}$ contains a
$K_{3}$ and a  $\overline {K_{3}}$ with some  vertex in common.
\end {claim}
 The proof is immediate and is omitted (in fact $U$ contains at least five 
$3$-element homogeneous subsets).

  \begin{claim} \label {syntaxic}A graph $U$ belongs to $Forb\{K_{1,3}, {\overline {K_{1,3}}}\}$ if and
only if 
$U(\{x, a\})=U(\{x, b\}) \not = U(\{y, a\})=U(\{y, b\})$ imply
$U(\{x, y\})=U(\{a, b\})$ for every $4$-tuple $(x, y, a, b)$ of distinct vertices. 
\end{claim} 
  \Proof  Immediate.\endproof   

In Claim \ref{voisinage} to Claim  \ref{degre4} below, we consider a graph $U\in Forb\{K_{1,3}, {\overline {K_{1,3}}}\}$.

\begin{claim}  \label {voisinage}The inequality  $\vert N_{U}(x)\cap N_{\overline
U}(y)\vert \leq 2$ holds for every unordered pair $\{x,y\}$ of distinct vertices.
\end{claim}
\Proof Let $A:= N_{U}(x)\cap N_{\overline U}(y)$. According to Claim  \ref {syntaxic},
$U(\{a, b\})=U(\{x,y\}) $ for every unordered pair
$\{a, b\}\in [A]^2$. Hence, if 
$\vert A\vert \geq 3$, $A$ contains either a $K_{3}$ or a $\overline {K_{3}}$. In particular $U$ 
contains either a 
$\overline {K_{1,3}}$ or a  $K_{1,3}$,  a contradiction.\endproof

\begin{claim} \label {bull1}  If $U$ contains a  $K_{3}$ and a $\overline {K_{3}}$ with a  vertex in
common then the graph induced by $U$ on the union of their vertices is a bull.
\end{claim}

\Proof
Let  $A:=\{x,a, b\}$ and   $B:= \{x, c, d\}$ be the vertex set of a $K_{3}$ and a
$\overline {K_{3}}$ respectively.    Since $U$ is claw-free then
$\overline U$ contains at least one edge from the vertex $a$ to the set $\{c, d\}$ and  an edge from 
the vertex $b$ to the set $\{c,d\}$. Two such edges cannot meet, otherwise $\overline U$
would not be claw-free. For the same reason $U$ contains two 
 disjoint edges from $\{c, d\}$ to $\{a,b\}$. The resulting
graph is a bull as claimed. 
\endproof          
 
 \begin{claim} \label{2k3k3bar} Let $x$ be a vertex belonging to a 
$\overline {K_{3}}$. Then  two $K_{3}$ containing
$x$ cannot have an edge in common.  \end{claim} 

 \Proof Supppose the  contrary.  Let $\{x,e, f\}$ be the vertices of a $\overline {K_{3}}$, 
$\{x,a, b\}$ and $\{x,b, c\}$ be the vertices of two
$K_{3}$. According to Claim  \ref{bull1} the graph induced by $U$ on $\{x, a, b, e, f\}$ is a
bull. With no loss of generality, we may suppose that 
$\{a, e\}$ and
$\{b,f\}$ are two edges of $U$, in which case $\{a, f\}$ and $\{b,
e\}$ are not edges of $U$. The graph induced by $U$ on 
$\{x, b, c, e, f\}$ is a bull too. Hence 
$\{c, e\}$ must be an  edge of $U$ whereas 
  $\{c, f\}$ is  not an edge of $U$. It follows that 
$\{x,a, c\}\subseteq N_{U}(b)\cap N_{\overline U}(f)$,
contradicting Claim \ref {voisinage}.
 \endproof
 
\begin{claim} \label{degre4} If a vertex $x$ belongs to a 
$\overline {K_{3}}$ then
$d_{U}(x)\leq 4$. Moreover,  the graph $W$ induced by $U$
on the union of the vertices of the $\overline {K_{3}}$ and the
neighbourhood of $x$ is isomorphic to $ \overline {A_{6}}$ or to  $\overline {E_{6}}$ if $d_{U}(x)=3$
and to $\overline {P_{9-\varepsilon}}$ if $d_{U}(x)=4$.
\end{claim}

\Proof
 Let $\{x,e, f\}$ be the vertices of a $\overline {K_{3}}$.

Let $U'$ be the graph
induced by $U$ on $N_{U}(x)$. According to Claim \ref{2k3k3bar}, $d_{U'}(y)\leq 1$ for every vertex $y\in N_{U}(x)$. Since in addition $U'$ cannot contain a $\overline {K_{3}}$, $\vert V(U')\vert \leq 4$.  In fact,  $U'$ consists of an edge
and an isolated vertex if $\vert V(U')\vert  =3$ and it consists of two disjoint edges
if $\vert V(U')\vert =4$. We have $V(U')=N_{U}(x)$, hence  $\vert V(U')\vert =d_U(x)$ and thus $d_U(x)\leq 4$.
Suppose $d_U(x)=3$. Let $\{a,b,c\}:=N_{U}(x)$; according to Claim \ref{bull1} 
we may suppose that the graph induced by $U$ on $\{x, a, b, e,f\}$ is a bull.  Since the graph induced on
$\{x,e,f\}$ is a $\overline {K_{3}}$, and $c$ is a  neighbour of $x$
then at  most one of the unordered pairs
$\{c,e\}$,  $\{c,f\}$ is an edge. If none
of these unordered pairs is an edge 
 then $W$ is isomorphic to $\overline
{A_{6}}$, otherwise $W$ is isomorphic to $\overline {E_{6}}$
as claimed. 
 Suppose $d_{U}(x)=4$. Let $\{a,b,c,d\}:=N_{U}(x)$;
we may suppose that the graph induced by $U$ on $\{x, a, b, e,
f\}$ is the bull  with $\{a,e\}, \{b,f\}\in E(U)$. The graph induced on $\{x, c, d, e,
f\}$ is a bull too. We may suppose  $\{e,c\}, \{d, f\}\in E(U)$
and $\{c ,f\}, \{d,e\}\not \in E(U)$.

In this case  $W$ is isomorphic to $\overline {P_{ 9-\varepsilon}}$ as
claimed. \endproof 
\subsubsection{Proof of Lemma \ref {caseanalysis}.}\label{Prooflemmacaseanalysis}Let $U\in Forb\{K_{1,3}, {\overline {K_{1,3}}}\}$ and
suppose that $U$ contains a
$K_{3}$ and a $\overline {K_{3}}$. Let
$V:=V(U)$.  
Then, obviously,   $\vert V\vert\geq 5$.  
Claim  \ref{disjointtriangle} insures that $U$ contains a $K_{3}$
and a $\overline {K_{3}}$ with a common vertex. Let $x$ be such a
vertex; according to Claim \ref {degre4}, we have
$d_{U}(x)\leq 4$ and by the same 	claim  applied to 
$\overline U$, we also have 
$d_{\overline U}(x)\leq 4$, proving
$\vert V\vert \leq 9$. \\

{\bf Case 1.}  $\vert V\vert=5$.  From Claim 
\ref{bull1},
$U$ is a bull. \\
 
{\bf Case 2.} $\vert V\vert = 6$. Then  either
$d_{U}(x)=3$ or
$d_{\overline U}(x)=3$. In the first case, $V$ is the union of the vertices of the $\overline {K_{3}}$ and the
neighbourhood of $x$, hence by Claim \ref {degre4}, 
$U$ is isomorphic to $ \overline {A_{6}}$ or to
$\overline {E_{6}}$.  In the second case, the same argument yields that $U$ is isomorphic  to
$A_{6}$ or to $E_{6}$.\\ 

{\bf Case 3.} $\vert V\vert = 7$.   Then  either (a)
$d_{U}(x)=4$ or
(b) $d_{\overline U}(x)=4$ or (c) $d_{U}(x)=d_{\overline U}(x)=3$.
According to Claim \ref {degre4}, in  case (a)
$U$ is isomorphic to $\overline {P_{9 -\varepsilon}}$ and in case
(b) to
$P_{9 -\varepsilon}$. Case (c) reduces to (a) or (b). Indeed, suppose that Case (c) holds. Let
$ \{a, b, c\}:=N_U(x)$ and 
$\{e, f, g\}:=N_{\overline U}(x)$.  We may suppose that $\{a, b\}\in E(U)$, $\{e, f\}\not \in E(U)$ and the
graph induced by $U$ on $\{x,a,b,e,f\}$ is the
bull considered in  the proof of Claim \ref{degre4}. Moreover, by Claim \ref {2k3k3bar},  
$\{c,a\},\{c,b\}\not \in E(U)$.  Since $\{c,x,e,f\}$ cannot  form a $K_{1,3}$, 
 either $\{c,e\}$ or
$\{c,f\}$ do not belong to $E(U)$. 
Suppose $\{c,f\}\not \in E(U)$ (the case $\{c,e\}\not
\in E(U)$ will be similar). Then $a$ belongs to a $K_{3}$  and a
$\overline {K_{3}}$. If $\{a,g\}\in E(U)$, then
$d_{U}(a)=4$ hence, from case (a), $U$ is
isomorphic to $\overline {P_{ 9 -\varepsilon}}$. If $\{a,g\}\not \in
E(U)$, then necessarily $\{b,
g\}\in E(U)$ (otherwise $\{x,a,b,g\}$ forms a
$\overline {K_{1,3}}$). 
 Necessarily $\{f,g\}\in E(U)$, otherwise $\{b,a,f,g\}$ forms a $K_{1,3}$. 
 Hence, $f$ belongs to a $K_{3}$
and a $\overline {K_{3}}$, also  $d_{\overline
U}(f)=4$. Thus  from case (b), $U$ is
isomorphic to $ P_{9-\varepsilon}$.\\

{\bf Case 4.} $\vert V\vert =8$. In this case, 
either 
  $d_{U}(x)=4$ or
$d_{\overline U}(x)=4$. Suppose 
$d_{U}(x)=4$ (otherwise, replace $U$ by  $\overline  
U$). Let $ \{a, b, c, d\}:=N_U(x)$ and 
$\{e, f, g\}:=N_{\overline U}(x)$. We may suppose that $\{a,b\}\in E(U)$, $\{e,f\}\not
\in E(U)$. According to Claim \ref {2k3k3bar},  $\{u,v\}\notin E(U)$ whenever $u\in \{a,b\}$ and 
  $v\in \{c,d\}$, thus $\{c,d\}\in E(U)$. 
 Let   $U_{-g}$ be the graph induced by $U$ on $V \setminus
\{g\}$. It satisfies (a) of  Case
3.  Hence it is  isomorphic to
$\overline {P_{9 -\varepsilon}}$ and we may suppose  that $E(U_{-g})=
A$, where $$A:=\{\{x,a\},\{x,b\},\{x,c\},\{x,d\},\{a,b\}, \{c,d\}, \{e,a\},\{e,c\},\{f,b\},\{f,d\}\}.$$
According to Claim \ref {2k3k3bar} applied to $\overline U$,  $\{g,e\}, \{g,f\}\in E(U)$. If $\{a,g\}\in E(U)$
also $\{d,g\}\in E(U)$ (otherwise, $\{a, g, e, d \}$ forms a $\overline {K_{1,3}}$),
$\{c,g\}\notin E(U)$ (otherwise, $\{ g,f,c,a \}$ forms a ${K_{1,3}}$),
$\{b,g\}\notin E(U)$ (otherwise, $\{ g,b,e,d \}$ forms a ${K_{1,3}}$), hence $U$ is isomorphic to 
$P_{9 -v}$.  If $\{a,g\}\notin E(U)$, then  $\{b,g\}\in E(U)$ (otherwise, $\{a, b, x, g \}$ forms a $\overline {K_{1,3}}$), 
then $\{d,g\}\notin E(U)$ (otherwise, $\{ d,g,b,e \}$ forms a ${K_{1,3}}$),
 $\{c,g\}\in E(U)$ (otherwise, $\{c,b,f, g \}$ forms a $\overline {K_{1,3}}$). Thus $U$ is isomorphic to  $P_{9 -v}$.\\

{\bf Case 5. }  $\vert V\vert =9$. In this case
$d_{U}(x)=4$. Let $\{a, b, c, d\}:=N_U(x)$ and 
$\{e, f, g, h\}:=N_{\overline U}(x)$.  Let  
$U_{-h}$ be the graph induced by $U$ on $V \setminus
\{h\}$. It satisfies the hypothesis of Case
4.  Hence it is  isomorphic to
${P_{9 -v}}$. Up to a relabeling, we may
suppose  that:
\begin{equation}\label{eq:paley}
E(U_{-h} )=A\cup \{\{g, a\}, \{g,d\},
\{g,e\}, \{g, f\}\}
\end{equation}

 Since $a$ belongs to a $\overline {K_{3}}$, then, by claim \ref{degre4}, 
 we have  
$d_{U}(a)\leq 4$, hence $\{a,h\}\not
\in E(U)$. Similarly, $\{d, h\}\not \in E(U)$ and $\{g, h\}\not \in E(U)$.  Hence, the unordered pairs $\{c,h\}, \{b,h\}, \{f,h\},  \{e,h\}$ belong to $E(U)$.  Taking account of Formula (\ref{eq:paley}), this yields that
$U$ is isomorphic to $P_{9}$. \endproof

\subsection{Ingredients for the proof of Theorem \ref{S(U)}.}
The proof of  the equivalence between Properties (1) and (2) of Theorem \ref{S(U)} relies on the
following lemma.

 \begin{lemma} \label{G,G',U}   Let $G$ and  $G'$ be two graphs on the same vertex set $V$ and let 
  $U:=G\dot{+}G'$. Then, the following properties are equivalent:
  \begin{enumerate}[{(a)}]
 \item     $G$ and $G'$ have the same $3$-element homogeneous subsets;
 \item $U(\{x,y\})=U(\{x,z\})\neq U(\{y,z\}) \Longrightarrow   G(\{x,y\})\neq G(\{x,z\})$ for all distinct elements $x,y,z$ of $V$. 
  \item  The sets $A_1:=E(U)\cap E(G)$ and  $A_2:=E(U)\setminus E(G)$  divide
$V(S({U}))$ into two
   independent sets and also 
    the sets $B_1:=E({\overline U})\cap E(G)$ and 
$B_2:=E({\overline U})\setminus E(G)$  divide
    $V(S({\overline U}))$ into two  independent sets.
 \end{enumerate}\end{lemma}

 \Proof   Observe first that Property (b) is equivalent to
the conjunction of the following properties:\\
$ (b_{U})$: If $\{u,v\}$ is an edge of $S(U)$ then $u\in E(G)$ iff      $v\notin E(G)$.            \\
 and\\
$  (b_{\overline U})$: If $\{u,v\}$ is an edge of $S(\overline U)$ then $u\in E(G)$ iff      $v\notin E(G)$.  

\noindent  $(a)  \Longrightarrow  (b)$. Let us show $(a)  \Longrightarrow (b_{U})$.\\
Let  $\{u,v\} \in E(S(U))$, then $u,v\in E(U)$.  
 By contradiction, we may suppose  that $u,v\in E(G)$ (the other case implies $u,v\in E(G')$ 
  thus is similar).   Since     $u$ and $v$ are edges
 of $U=G\dot{+}G'$ then   $u,v\notin E(G')$. Let $w:=\{y,z\}$ such that   $u=\{x,y\}$,  $v=\{x,z\}$. Then
 $w\notin E(U)$ and thus $w \in E(G)$ iff   $w\in E(G')$. \\
 If $w \in E(G)$,
 $\{x,y,z\}$  is  an homogeneous subset of $G$. Since $G$ and $G'$ have the same $3$-element homogeneous subsets, 
  $\{x,y,z\}$  is an homogeneous subset of $G'$. Hence, since $u,v\notin E(G')$, 
    $w=\{y,z\}     \notin E(G')$, thus $w\notin E(G)$, a contradiction.\\
If $w  \notin E(G)$, then $w  \notin E(G')$; since $u,v \notin E(G')$ it follows that $\{x,y,z\}$  is an homogeneous subset of $G'$. Consequently $\{x,y,z\}$  is  an homogeneous subset of $G$.
 Since $u,v\in E(G)$, then $w  \in E(G)$, a contradiction.\\
The implication    $a)  \Longrightarrow (b_{\overline U})$ is similar.\\
$ (b)  \Longrightarrow (a)$. Let $T$ be a $K_{3}$ of $G$.
Suppose that $T$ is not an homogeneous subset of $G'$ then we may suppose 
$T=\{u,v,w\}$ with $u,v\in E(G')$ and $w\notin E(G')$
or   $u,v\in E(\overline {G'})$ and
  $w\notin E(\overline {G'})$. In the first case $\{u,v\}\in E(S(\overline U))$, which contradicts Property
   $ (b_{\overline U})$,
  in the second case     $\{u,v\} \in E(S(U))$, which contradicts  Property $ (b_{U})$.\\
 $ (b)    \Longrightarrow (c)$. First $V(S(U)) =E(U)=A_1\cup A_2$ and
 $V(S(\overline U))=E(\overline U)=B_1\cup B_2$.
 Let  $u,v$ be two distinct elements 
of $A_1$ (respectively  $A_2$). Then  $u,v\in E(G)$ (respectively   $u,v\notin E(G)$). From $ (b_{U})$  we have 
 $\{u,v\}\notin E(S(U))$. Then $A_1$ and $A_2$ are   independent sets of $V(S(U))$.
 The proof  that $B_1$ and $B_2$ are  independent sets of $V(S({\overline U}))$ is similar.
 
$(c)  \Longrightarrow (b)$. This implication is trivial. \endproof
\subsection {Proof of Theorem \ref{S(U)}.}\label{Proof of TheoremS(U)}
Implication $ (1)  \Longrightarrow (2)$ follows directly from  implication $(a)  \Longrightarrow (c)$ of Lemma \ref{G,G',U}.  Indeed,  Property (c) implies trivially that  $S(U)$ and $S(\overline U)$ are bipartite. \\
$ (2)  \Longrightarrow (1)$.     Suppose   that $S(U)$ and  $S(\overline U)$  are bipartite.   Let $\{A_1, A_2\}$ and
$\{B_1, B_2\}$  be respectively a partition of   $V(S(U))=E(U)$ and  $V(S(\overline U))=E(\overline U)$ into independent sets. Note that $A_i\cap B_j =\emptyset$, for $i,j\in \{1,2\}$. 
Let $G,G'$ be two graphs with the same vertex set as $U$ such that $E(G)=A_1\cup B_1$ and     $E(G')=A_2\cup B_1$. 
Clearly  $E(G\dot{+}G')=A_1\cup A_2=E(U)$. Thus $U=G\dot{+}G'$. 
To conclude that Property $(1)$ holds, it suffices to  show that $G$ and $G'$ have the same $3$-element homogeneous subsets, that is Property (a) of Lemma \ref{G,G',U} holds. For that, note that $A_1=E(U)\cap E(G)$, $A_2=E(U)\setminus E(G)$, $B_1=E(\overline U)\cap E(G)$ and $B_2=E(\overline U)\setminus E(G)$ and thus Property (c) of 
Lemma \ref{G,G',U} holds. It follows that Property (a) of this lemma holds.

The proof of implication $(2)  \Longrightarrow (3)$ was given in Section 1. For the converse implication, let $U$ be a graph satisfying Property (3). It is clear from  Figure 1  that  $S(P{_9})$ is bipartite. Since $\overline {P_{9}}$ is isomorphic to $P_9$,  $S(\overline {P{_9}})$ is bipartite too. Thus, if $U$ is isomorphic to an induced subgraph of $P_{9}$, Property (2) holds. If not, we may suppose that  the connected components of $U$ are cycles of even length, paths or isolated vertices (otherwise, replace $U$ by $\overline U$). In this case, $S(U)$ is trivially bipartite. In order to prove that Property 2 holds, it suffices to prove that  $S(\overline U)$ is bipartite too. This is a direct consequence of the following claim:

\begin{claim} If $U$ is a bipartite graph, then $S(\overline U)$ is bipartite too.
\end{claim}

\Proof  If $c: V(U) \rightarrow \Z/2\Z$ is a colouring of $U$, set $c':
V(S(\overline U))\rightarrow \Z/2\Z$ defined by $c'(\{x, y\}):= c(x)+c(y)$. \endproof

 With this, the proof of Theorem \ref{S(U)} is complete. \endproof \\

\noindent{\bf Acknowledgements}\\

We  thank S. Thomass\'e for his helpful
comments. We thank the referees for their suggestions.

\end{document}